\documentclass[12pt,a4paper]{article}
\usepackage[dvips]{graphicx}

\usepackage{amssymb}
\usepackage{amsmath}
\usepackage{amsfonts}
\usepackage{amsthm,amscd}
\usepackage{amsbsy}
\usepackage{latexsym}

\theoremstyle{plain}
\newtheorem{Thm}{Theorem}[section]

\newtheorem{Prop}[Thm]{Proposition}
\theoremstyle{definition}
\newtheorem{Def}[Thm]{Definition}
\newtheorem{Rem}[Thm]{Remark}
\newtheorem{Ex}[Thm]{Example}

\newcommand{\Proof}[2][Proof]{\begin{proof}[{#1}] #2 \end{proof}}

\newcommand{\Supp}{\mathop{\rm Supp}}

\newcommand{\eps}{\ensuremath{\varepsilon}}

\renewcommand{\tilde}{\widetilde}
\renewcommand{\bar}{\overline}

\newcommand{\rbra}[1]{\!\left( #1 \right)} 
\newcommand{\cbra}[1]{\!\left\{ #1 \right\}} 

\newcommand{\bmat}[1]{\begin{bmatrix} #1 \end{bmatrix}} 

\newcommand{\bR}{\ensuremath{\mathbb{R}}}
\newcommand{\bZ}{\ensuremath{\mathbb{Z}}}

\newcommand{\cP}{\ensuremath{\mathcal{P}}}

\newcounter{Section}
\newcommand{\seclabel}[1]{\refstepcounter{Section} \label{#1}}

\begin{document}

\begin{center}
{\Large \bf Random walk in a finite directed graph subject to a synchronizing road coloring 
}
\end{center}

\begin{center}
Kouji Yano\footnote{
Graduate School of Science, Kyoto University, Kyoto, JAPAN.}\footnote{
The research of this author was supported by KAKENHI (20740060).} 
\quad and \quad 
Kenji Yasutomi\footnote{
Graduate School of Science and Engineering, 
Ritsumeikan University, Kusatsu, JAPAN.} 
\end{center}


\begin{abstract}
A constructive proof is given to the fact that 
any ergodic Markov chain can be realized as a random walk 
subject to a synchronizing road coloring. 
Redundancy (ratio of extra entropy) in such a realization is also studied. 
\end{abstract}

\section{Introduction} \seclabel{sec: intro}

A random walk in $ \bR $ 
is a process $ (S_n)_{n \ge 0} $ 
which may be represented as 
\begin{align}
S_n = \xi_n + \xi_{n-1} + \cdots + \xi_1 + S_0 
, \quad n \ge 1 
\label{eq: Sn}
\end{align}
for some sequence $ (\xi_n)_{n \ge 1} $ of IID 
(i.e., independent and identically distributed) random variables 
being independent of $ S_0 $. 
Note that equation \eqref{eq: Sn} is equivalent to the recursion relation 
\begin{align}
S_n = \xi_n + S_{n-1} 
, \quad n \ge 1 . 
\label{}
\end{align}

We may introduce a natural analogue of random walk 
taking values in a finite set $ V $, say, $ \{ 1,\ldots,m \} $. 
Let $ \Sigma $ denote the set of all mappings of $ V $ into itself. 
A {\em random walk in $ V $} 
is a pair of processes $ \{ (X_n)_{n \ge 0},(\phi_n)_{n \ge 1} \} $ 
such that 
$ (\phi_n)_{n \ge 1} $ is a sequence of IID random variables taking values in $ \Sigma $ 
and being independent of $ X_0 $ and such that 
\begin{align}
X_n = (\phi_n \circ \phi_{n-1} \circ \cdots \circ \phi_1)(X_0) 
, \quad n \ge 1 . 
\label{eq: Xn}
\end{align}
Note that equation \eqref{eq: Xn} is equivalent to the recursion relation 
\begin{align}
X_n = \phi_n(X_{n-1}) 
, \quad n \ge 1 . 
\label{}
\end{align}
It is obvious that, for each $ n \ge 1 $, 
the random variable $ \phi_n $ is independent of $ \sigma(X_j,\phi_j:j \le n-1) $, 
since each $ X_j $ is measurable with respect to 
$ \sigma(X_0,\phi_1,\ldots,\phi_j) $. 

It is now natural to extend the index set to $ \bZ $, the set of all integers, as follows. 

\begin{Def}
A {\em random walk in $ V $ parametrized by $ \bZ $} 
is a pair of processes $ \{ (X_n)_{n \in \bZ},(\phi_n)_{n \in \bZ} \} $ 
which satisfies the following conditions: 
\begin{enumerate}
\item 
$ (\phi_n)_{n \in \bZ} $ is a sequence of IID random variables taking values in $ \Sigma $; 
\item 
for each $ n \in \bZ $, the random variable $ \phi_n $ 
is independent of $ \sigma(X_j,\phi_j:j \le n-1) $; 
\item 
$ X_n = \phi_n(X_{n-1}) $ holds almost surely for all $ n \in \bZ $. 
\end{enumerate}
If $ \phi_n $'s have common law $ \mu $ on $ \Sigma $, 
such a random walk is called a {\em $ \mu $-random walk}. 
\end{Def}

Our $ \mu $-random walk may also be called 
a {\em random walk in a finite directed graph subject to a road coloring}. 
The reason will be explained in Section \ref{sec: rc}. 
Each element of $ V $ will be called a {\em site}. 

If a $ \mu $-random walk $ \{ (X_n)_{n \in \bZ},(\phi_n)_{n \in \bZ} \} $ is given, 
then the process $ (X_n)_{n \in \bZ} $ is a Markov chain in $ V $ 
whose one-step transition probability is given by 
\begin{align}
P(X_1=y|X_0=x) = \mu \rbra{ \sigma \in \Sigma : \sigma(x)=y } 
, \quad x,y \in V . 
\label{}
\end{align}
Conversely, if a Markov chain $ (Y_n)_{n \in \bZ} $ is given, 
we call $ \mu $ a {\em mapping law} for the Markov chain 
if the following identity holds: 
\begin{align}
P(Y_1=y|Y_0=x) = \mu \rbra{ \sigma \in \Sigma : \sigma(x)=y } 
, \quad x,y \in V . 
\label{}
\end{align}
In this case, for any $ \mu $-random walk $ \{ (X_n)_{n \in \bZ},(\phi_n)_{n \in \bZ} \} $, 
we can easily show that 
the Markov chain $ (X_n)_{n \in \bZ} $ is identical in law to $ (Y_n)_{n \in \bZ} $. 

\begin{Prop} \label{prop: mplaw}
For any Markov chain $ (Y_n)_{n \in \bZ} $ in $ V $, there exists a mapping law $ \mu $. 
\end{Prop}

A proof by means of rational approximation can be found in \cite{YanoYasu}. 
We shall give its constructive proof in the next section. 

Note that, for a given Markov chain, there may exist several mapping laws. 
We may expect that we can take a nice mapping law in the following sense. 

\begin{Def} \label{def: sync mu}
Let $ \mu $ be a probability law on $ \Sigma $ 
and denote by $ \Supp(\mu) $ the support of $ \mu $. 
We say that $ \mu $ is {\em synchronizing} (or simply {\em sync}) if 
there exists a finite sequence $ \sigma_1,\ldots,\sigma_p $ 
of elements of $ \Supp(\mu) $ such that 
$ \sigma_p \circ \sigma_{p-1} \circ \cdots \circ \sigma_1 $ maps $ V $ into a signleton. 
\end{Def}

Note that a $ \mu $-random walk associated with a sync mapping law 
is utilized in Propp--Wilson's sampling method of stationary law, 
which is called {\em coupling from the past}; 
we shall mention it briefly in Section \ref{sec: cftp}. 

Suppose that $ \mu $ is sync 
and let $ \{ (X_n)_{n \in \bZ},(\phi_n)_{n \in \bZ} \} $ be a $ \mu $-random walk. 
We may assume without loss of generality that 
for any $ x \in V $ there exists $ \sigma \in \Supp(\mu) $ such that $ x \in \sigma(V) $; 
in fact, the Markov chain $ (X_n)_{n \in \bZ} $ never visits such sites $ x \in V $ 
that $ x \notin \sigma(V) $ for any $ \sigma \in \Supp(\mu) $. 
Then we see that the Markov chain $ (X_n)_{n \in \bZ} $ is {\em ergodic}, 
i.e., the following two conditions hold (see, e.g., \cite{Nor}): 
\begin{enumerate}
\item the Markov chain is {\em irreducible}, i.e., 
$ P(X_0=x)>0 $ for all $ x \in V $ 
and for any $ x,y \in V $ there exists $ n \ge 1 $ such that $ P(X_n=y|X_0=x)>0 $; 
\item the Markov chain is {\em aperiodic}, i.e., 
for any $ x \in V $, 
the greatest common divisor of $ \{ n \ge 1:P(X_n=x|X_0=x)>0 \} $ is one. 
\end{enumerate}
The condition (i) is obvious. The condition (ii) may be verified as follows. 
Let $ x \in V $. 
Take $ a \in V $ such that $ \sigma_p \circ \cdots \circ \sigma_1(V) = \{ a \} $ 
and take $ q \ge 1 $ such that $ P(X_q=x|X_0=a)>0 $. 
Then the set $ \{ n \ge 1:P(X_n=x|X_0=x)>0 \} $ contains 
all integers greater than $ p+q $, and hence its greatest common divisor is one. 

The following theorem asserts that the converse is also true. 

\begin{Thm}[\cite{YanoYasu}] \label{thm: sync}
Suppose that $ (Y_n)_{n \in \bZ} $ is an ergodic Markov chain. 
Then there exists a sync mapping law. 
\end{Thm}

To prove Theorem \ref{thm: sync}, the authors in \cite{YanoYasu} 
utilized a profound graph-theoretic theorem, 
which was recently obtained by Trahtman \cite{MR2534238}, 
the complete solution to the {\em road coloring problem}; 
we shall explain it briefly in Section \ref{sec: rc}. 
In this chapter, we would like to give an elementary, self-contained and constructive proof 
of Theorem \ref{thm: sync} 
without using Trahtman's theorem. 

The remainder of this chapter is as follows. 
In Section \ref{sec: prf}, we give constructive proofs 
to Proposition \ref{prop: mplaw} and Theorem \ref{thm: sync}. 
In Section \ref{sec: red}, we study redundancy in random walk realization of a Markov chain. 
In Section \ref{sec: cftp}, we mention the couplinng from the past. 
In Section \ref{sec: rc}, we explain how our random walk is related to road coloring. 
In Section \ref{sec: conc}, we provide a summary and conclusion.

\section{A constructive proof of existence of sync mapping law} \seclabel{sec: prf}

A matrix $ Q=(q_{x,y})_{x,y \in V} $ with non-negative entries 
is called a {\em transition matrix} if 
\begin{align}
\sum_{y \in V} q_{x,y} = 1 
\quad \text{for all $ x \in V $}. 
\label{eq: trans mat}
\end{align}
We give a constructive proof of Proposition \ref{prop: mplaw} for later use. 

\Proof[A constructive proof of Proposition \ref{prop: mplaw}]{
It suffices to show that, for any transition matrix $ Q $, 
there exists a mapping law $ \mu $ for $ Q $, i.e., 
\begin{align}
q_{x,y} = \mu(\sigma \in \Sigma : \sigma(x)=y) 
, \quad x,y \in V . 
\label{eq: mplaw}
\end{align}
We define 
\begin{align}
E(Q) = \{ (x,y) \in V \times V : q_{x,y}>0 \} . 
\label{eq: EQ}
\end{align}
Let us prove the result by induction of $ \sharp E(Q) $, 
where $ \sharp A $ stands for the number of elements of $ A $. 
It is obvious by \eqref{eq: trans mat} that 
$ \sharp \{ y \in V : q_{x,y}>0 \} \ge 1 $ for all $ x \in V $, 
and hence that $ \sharp E(Q) \ge \sharp V $. 

Suppose that $ \sharp E(Q)=\sharp V $. 
Then, by \eqref{eq: trans mat}, it holds that 
$ \sharp \{ y \in V : q_{x,y}>0 \} = 1 $ for all $ x \in V $. 
This shows that there exists $ \sigma \in \Sigma $ such that 
\begin{align}
q_{x,y} = 
\begin{cases}
1 & \text{if $ y=\sigma(x) $}, \\
0 & \text{otherwise}. 
\end{cases}
\label{}
\end{align}
Thus the Dirac mass at $ \sigma $ is as desired. 

Let $ k > \sharp V $ and suppose that 
all transition matrix $ Q $ such that $ \sharp E(Q) < k $ admits a mapping law. 
Let $ Q $ be a transition matrix such that $ \sharp E(Q) = k $. 
Write 
\begin{align}
\eps = \min \{ q_{x,y} : (x,y) \in E(Q) \} . 
\label{}
\end{align}
Since $ \sharp E(Q)>\sharp V $, we see that $ 0<\eps<1 $. 
Take $ (x,y) \in E(Q) $ such that $ \eps = q_{x,y} $, 
and take $ \sigma \in \Sigma $ such that $ (x,\sigma(x)) \in E(Q) $ for all $ x \in V $. 
Define $ \tilde{Q}=(\tilde{q}_{x,y})_{x,y \in V} $ by 
\begin{align}
\tilde{q}_{x,y} = \frac{1}{1-\eps} \rbra{ q_{x,y} - \eps 1_{\{ \sigma(x)=y \}} } . 
\label{}
\end{align}
Then we see that $ \tilde{Q} $ is a transition matrix 
and that $ \sharp E(\tilde{Q}) < k $. 
Now, by the assumption of the induction, we see that 
$ \tilde{Q} $ admits a mapping law $ \tilde{\mu} $. 
Therefore we conclude that $ (1-\eps) \tilde{\mu} + \eps \delta_{\sigma} $ 
is a mapping law for $ Q $. 
The proof is now complete. 
}

Utilizing Proposition \ref{prop: mplaw}, 
we give a constructive proof of Theorem \ref{thm: sync}.

\Proof[A constructive proof of Theorem \ref{thm: sync}]{
Let $ Q $ be the transition matrix of an ergodic Markov chain. 
Then we see that there exists $ r \ge 1 $ such that 
the $ r $-th product $ Q^r $ has positive entries. 

Take $ x_0 \in V $ arbitrarily and set $ V_r = \{ x_0 \} $. 
If $ V_k $ is defined for $ k=r,r-1,\ldots,1 $, define $ V_{k-1} $ recursively by 
\begin{align}
V_{k-1} = \cbra{ x \in V : (x,y) \in E(Q) \ \text{for some $ y \in V_k $} } , 
\label{}
\end{align}
where $ E(Q) $ has been defined in \eqref{eq: EQ}. 
Note that 
\begin{align}
\sharp V_r \le \sharp V_{r-1} \le \cdots \le \sharp V_0 . 
\label{}
\end{align}
Since $ Q^r $ has positive entries, we see that $ V_0 = V $. 

For $ k=r,r-1,\ldots,1 $, we pick $ \sigma_k \in \Sigma $ so that 
$ \sigma_k(x) \in V_k $ if $ x \in V_{k-1} $ 
and $ (x,\sigma_k(x)) \in E(Q) $ if $ x \notin V_{k-1} $. 
We then have 
\begin{align}
\sigma_{r} \circ \sigma_{r-1} \circ \cdots \circ \sigma_1(V) = \{ x_0 \} . 
\label{}
\end{align}
Let $ \mu^{(1)} $ denote the uniform law on $ \{ \sigma_1,\ldots,\sigma_{r} \} $. 
Then we see that $ \mu^{(1)} $ is a sync mapping law 
for a transition matrix $ Q^{(1)} = (q^{(1)}_{x,y})_{x,y \in V} $ where 
\begin{align}
q^{(1)}_{x,y} = \mu^{(1)}( \sigma \in \Sigma : y=\sigma(x) ) 
= \frac{1}{r} \sum_{k=1}^{r} 1_{\{ y=\sigma_k(x) \}} . 
\label{}
\end{align}
Write 
\begin{align}
\eps = \min \{ q_{x,y} : (x,y) \in E(Q) \} > 0 
\label{}
\end{align}
and define $ Q^{(2)}=(q^{(2)}_{x,y})_{x,y \in V} $ by 
\begin{align}
q^{(2)}_{x,y} = \frac{1}{1-\eps} \rbra{ q_{x,y} - \eps q^{(1)}_{x,y} } . 
\label{}
\end{align}
Then $ Q^{(2)} $ is a transition matrix, so that 
we may obtain a mapping law $ \mu^{(2)} $ for $ Q^{(2)} $ 
in the constructive way of the proof of Proposition \ref{prop: mplaw} given above. 

Now we define 
\begin{align}
\mu = \eps \mu^{(1)} + (1-\eps) \mu^{(2)} , 
\label{}
\end{align}
which we have proved that is a sync mapping law for $ Q $. 
The proof is therefore complete. 
}

\section{Redundancy in random walk realization} \seclabel{sec: red}

The uncertainty associated with information source may be measured by entropy 
(see, e.g., \cite{Bil}). 
A Markov chain $ Y=(Y_n)_{n \in \bZ} $ 
with transition matrix $ Q=(q_{x,y})_{x,y \in V} $ 
and with stationary law $ \lambda $ has its entropy given by 
\begin{align}
h(Y) = - \sum_{x,y \in V} \lambda(x) q_{x,y} \log q_{x,y} , 
\label{}
\end{align}
where we adopt the binary logarithm $ \log = \log_2 $ for simplicity, 
and follow the usual convention: $ 0 \log 0 = 0 $. 
For a probability law $ \mu $ on $ \Sigma $, 
an IID sequence $ \phi=(\phi_n)_{n \in \bZ} $ with common law $ \mu $ 
has its entropy given by 
\begin{align}
h(\phi) = h(\mu) = - \sum_{\sigma \in \Sigma} \mu(\sigma) \log \mu(\sigma) . 
\label{eq: hmu}
\end{align}
A $ \mu $-random walk $ (X,\phi) = \{ (X_n)_{n \in \bZ},(\phi_n)_{n \in \bZ} \} $ 
with stationary law $ \lambda $ 
is a Markov chain whose transition matrix 
$ \bar{Q}=(\bar{q}_{(x,\nu),(y,\sigma)})_{(x,\nu),(y,\sigma) \in V \times \Sigma} $ 
and stationary law $ \bar{\lambda} $ given by 
\begin{align}
\bar{q}_{(x,\nu),(y,\sigma)} 
= \mu(\sigma) 1_{\{ y=\sigma(x) \}} 
, \quad 
\bar{\lambda}((x,\nu)) = \mu(\nu) \lambda(w \in V : x=\nu(w)) . 
\label{}
\end{align}
Now its entropy $ h(X,\phi) $ is computed as 
\begin{align}
h(X,\phi) 
=& - \sum_{(x,\nu),(y,\sigma) \in V \times \Sigma} 
\bar{\lambda}((x,\nu)) \bar{q}_{(x,\nu),(y,\sigma)} \log \bar{q}_{(x,\nu),(y,\sigma)} 
\label{} \\
=& - \sum_{x,y \in V , \ \sigma \in \Sigma} 
\cbra{ \sum_{\nu \in \Sigma} \mu(\nu) \lambda(w \in V : x=\nu(w)) } 
1_{\{ y=\sigma(x) \}} 
\mu(\sigma) \log \mu(\sigma) 
\label{} \\
=& - \sum_{x \in V , \ \sigma \in \Sigma} 
\lambda(x) 
\cbra{ \sum_{y \in V} 1_{\{ y=\sigma(x) \}} } 
\mu(\sigma) \log \mu(\sigma) 
\label{} \\
=& - \cbra{ \sum_{x \in V} \lambda(x) } 
\cbra{ \sum_{\sigma \in \Sigma} \mu(\sigma) \log \mu(\sigma) } 
\label{} \\
=& - \sum_{\sigma \in \Sigma} \mu(\sigma) \log \mu(\sigma) . 
\label{}
\end{align}
Thus we obtain $ h(X,\phi)=h(\phi)=h(\mu) $. 

If the Markov chain $ Y $ is identical in law to $ X $ 
for some $ \mu $-random walk $ (X,\phi) $, we have 
\begin{align}
h(\mu) \ge h(Y) . 
\label{eq: ineq}
\end{align}
In fact, by \eqref{eq: mplaw}, 
we have $ \mu(\sigma) \le q_{x,y} $ if $ y=\sigma(x) $, 
and hence we see that 
\begin{align}
h(\mu) 
=& - \sum_{x \in V} 
\lambda(x) \sum_{\sigma \in \Sigma} \mu(\sigma) \log \mu(\sigma) 
\label{} \\
=& - \sum_{x,y \in V} 
\lambda(x) \sum_{\sigma \in \Sigma: \, y=\sigma(x)} \mu(\sigma) \log \mu(\sigma) 
\label{} \\
\ge& - \sum_{x,y \in V} 
\lambda(x) \sum_{\sigma \in \Sigma: \, y=\sigma(x)} \mu(\sigma) \log q_{x,y} 
\label{} \\
= & - \sum_{x,y \in V} 
\lambda(x) q_{x,y} \log q_{x,y} = h(Y) . 
\label{}
\end{align}
The inequality \eqref{eq: ineq} shows that any $ \mu $-random walk realization of $ Y $ 
requires some extra entropy, the extent of which may be measured by 
\begin{align}
r(\mu;Y) := \frac{h(\mu) - h(Y)}{h(\mu)} . 
\label{}
\end{align}
This quantity $ r(\mu;Y) $ is called the {\em (relative) redundancy} 
in the $ \mu $-random walk realization of the Markov chain $ Y $. 
We denote the totality of all possible redundancies by 
\begin{align}
\rho(Y) = \{ r(\mu;Y) : \ \text{$ \mu $ is a mapping law for $ Y $} \} . 
\label{}
\end{align}

\begin{Thm} \label{thm: rdd}
For a Markov chain $ Y $, the following assertions hold: 
\begin{enumerate}
\item 
the set $ \rho(Y) $ has finite minimum $ r(Y) \ge 0 $ and maximum $ R(Y) \le 1 $; 
\item 
for any $ r(Y) \le r \le R(Y) $, there exists a mapping law $ \mu $ for $ Y $ 
such that $ r(\mu;Y)=r $. 
\end{enumerate}
Moreover, if $ Y $ is ergodic, then the following assertion also holds: 
\begin{enumerate} \setcounter{enumi}{2}
\item 
for any $ r(Y) < r < R(Y) $, there exists a sync mapping law $ \mu $ for $ Y $ 
such that $ r(\mu;Y)=r $. 
\end{enumerate}
\end{Thm}

\Proof{
Let us remark on several basic facts about the entropy. 
Since $ \Sigma $ is a finite set, the totality of probability measures on $ \Sigma $, 
which is denoted by $ \cP(\Sigma) $, is equipped with the total variation topology. 
It is well-known that $ \cP(\Sigma) $ is compact and that 
$ \mu_n \to \mu $ if and only if $ \mu_n(\sigma) \to \mu(\sigma) $ for all $ \sigma \in \Sigma $. 
By definition \eqref{eq: hmu}, 
the function $ \cP(\Sigma) \ni \mu \mapsto h(\mu) $ is continuous. 

(i) Let $ \cP(Y) $ denote the set of all mapping laws for $ Y $. 
It is obvious that $ \cP(Y) $ is a compact convex subset of $ \cP(\Sigma) $. 
Since $ h(\mu) \ge h(Y) > 0 $ for all $ \mu \in \cP(Y) $, 
and since $ t \mapsto (t-h(Y))/t $ is continuous in $ t \ge h(Y) $, 
we see that $ \cP(Y) \ni \mu \mapsto r(\mu;Y) $ is continuous. 
Hence we see that the set $ \rho(Y) $ has finite minimum $ r(Y) $ and maximum $ R(Y) $. 

(ii) 
Take $ \mu^{(1)},\mu^{(2)} \in \cP(Y) $ such that 
$ r(Y)=r(\mu^{(1)};Y) $ and $ R(Y)=r(\mu^{(2)};Y) $. 
Let $ 0 \le p \le 1 $. 
Then $ \mu_p := p \mu^{(1)} + (1-p) \mu^{(2)} $ also belongs to $ \cP(Y) $. 
Since $ [0,1] \ni p \mapsto r(\mu_p;Y) $ is continuous, 
we see that $ \rho(Y) $ contains all $ r $ such that $ r(Y) < r < R(Y) $. 
Thus we obtain (ii). 

(iii) 
Suppose that $ Y $ is ergodic. 
Theorem \ref{thm: sync} implies that there exists a sync mapping law $ \mu^{(0)} $ for $ Y $. 
Let $ r(Y)<r<R(Y) $ and take $ r^{(1)} $, $ r^{(2)} $ such that 
$ r(Y)<r^{(1)}<r $ and $ r<r^{(2)}<R(Y) $. 
By (ii), we may take mapping laws $ \mu^{(1)} $ and $ \mu^{(2)} $ for $ Y $ 
such that $ r(\mu^{(1)};Y)=r^{(1)} $ and $ r(\mu^{(2)};Y)=r^{(2)} $. 
Now we may take $ \eps>0 $ small enough such that 
\begin{align}
r((1-\eps)\mu^{(1)}+\eps \mu^{(0)};Y) < r < r((1-\eps)\mu^{(2)}+\eps \mu^{(0)};Y) . 
\label{}
\end{align}
Hence we may take $ 0<p<1 $ such that the mapping law $ \mu $ defined by 
\begin{align}
\mu = (1-\eps)(p\mu^{(1)}+(1-p)\mu^{(2)}) + \eps \mu^{(0)} 
\label{}
\end{align}
satisfies $ r(\mu;Y)=r $. 
This shows that $ \mu $ is a sync mapping law for $ Y $. 
Therefore the proof is complete. 
}

\begin{Ex} \label{ex: example}
Let $ V = \{ 1,2,3 \} $ and let 
\begin{align}
Q = 
\bmat{
q_{1,1} & q_{1,2} & q_{1,3} \\
q_{2,1} & q_{2,2} & q_{2,3} \\
q_{3,1} & q_{3,2} & q_{3,3} 
}
= 
\bmat{
0   & 0   & 1   \\
1/2 & 0   & 1/2 \\
1/2 & 1/2 & 0 
}
. 
\label{}
\end{align}
The Markov chain $ Y $ with transition matrix $ Q $ 
has a unique stationary law 
\begin{align}
\lambda 
= \bmat{\lambda(1), \lambda(2), \lambda(3) } 
= \frac{1}{9} \bmat{3, 2, 4} . 
\label{}
\end{align}
A simple computation leads to $ h(Y) = 2/3 $. 
For a mapping law $ \mu $ for $ Y $, 
elements which may possibly be contained in $ \Supp(\mu) $ are the following four: 
\begin{align}
\sigma^{(1)} = 
\bmat{
0 & 0 & 1 \\
1 & 0 & 0 \\
1 & 0 & 0 
}, \quad 
\sigma^{(2)} = 
\bmat{
0 & 0 & 1 \\
0 & 0 & 1 \\
0 & 1 & 0 
}, \quad 
\sigma^{(3)} = 
\bmat{
0 & 0 & 1 \\
0 & 0 & 1 \\
1 & 0 & 0 
}, \quad 
\sigma^{(4)} = 
\bmat{
0 & 0 & 1 \\
1 & 0 & 0 \\
0 & 1 & 0 
}. 
\label{}
\end{align}
Set $ p = \mu(\sigma^{(1)}) $. A simple computation leads to 
\begin{align}
\mu(\sigma^{(2)}) = p 
, \quad 
\mu(\sigma^{(3)}) = 
\mu(\sigma^{(4)}) = 1/2 - p . 
\label{}
\end{align}
Thus we obtain 
\begin{align}
h(\mu) = 2f(p) + 2f(1/2-p) , 
\label{}
\end{align}
where $ f(t) = - t \log t $. 
Since the variable $ p $ may vary in $ [0,1/2] $, we see that 
$ h(\mu) $ ranges $ [1,2] $, 
where the minimum $ h(\mu)=1 $ is attained at $ p=0 $ and $ 1/2 $ 
and the maximum $ h(\mu)=2 $ at $ p=1/4 $. 
Hence we obtain 
\begin{align}
r(Y) = 1/3 , \quad R(Y) = 2/3 . 
\label{}
\end{align}
In this case, for all $ p \in [0,1/2] $, 
the mapping law $ \mu $ is sync; 
in fact, 
$ \sigma^{(1)} \circ \sigma^{(2)}(V) = \{ 1 \} $ 
and 
$ \sigma^{(3)} \circ \sigma^{(4)} \circ \sigma^{(4)} \circ \sigma^{(3)}(V) = \{ 3 \} $. 
\end{Ex}

Before closing this section, we mention the following theorem, 
which provides a necessary and sufficient condition for zero minimum redundancy. 

\begin{Thm}[\cite{YanoYasu}] \label{thm: punif}
Suppose that $ Y $ is ergodic. Then $ r(Y)=0 $ if and only if 
$ Y $ is {\em p-uniform}, i.e., 
there exists a probability law $ \nu $ on $ V $ 
and a family $ \{ \tau_x:x \in V \} $ of permutations of $ V $ such that 
\begin{align}
q_{x,y} = \nu(\tau_x(y)) , \quad x,y \in V . 
\label{}
\end{align}
\end{Thm}

For the proof of Theorem \ref{thm: punif}, see \cite{YanoYasu}.

\section{Coupling from the past} \seclabel{sec: cftp} 

In some practical problems, 
we sometimes need to simulate the stationary law of an ergodic Markov chain. 
As a powerful method for the simulation, 
Propp--Wilson's coupling from the past is widely known; 
see \cite{MR1611693} and also \cite{Hag} and \cite{LPW}. 
The fundamental idea is to utilize 
a random walk realization associated with a sync mapping law. 
Let us explain it briefly. 

Let an ergodic Markov chain be given 
and suppose that we find a sync mapping law $ \mu $ for the Markov chain. 
Then a $ \mu $-random walk $ \{ (X_n)_{n \in \bZ},(\phi_n)_{n \in \bZ} \} $ 
is a realization of the Markov chain. 
Let $ (\sigma_1,\ldots,\sigma_p) $ be a finite sequence of elements of $ \Supp(\mu) $ such that 
$ \sigma_p \circ \cdots \circ \sigma_1(V) $ is a singleton. 
The latest time 
when the exact sequence $ (\sigma_p,\sigma_{p-1},\ldots,\sigma_1) $ 
can be found in $ (\phi_0,\phi_{-1},\ldots) $ will be denoted by 
\begin{align}
T = \sup \{ k \in \bZ : 0 \ge k+p-1 , \ \phi_{k+p-1}=\sigma_p, \ldots, \phi_k=\sigma_1 \} . 
\label{}
\end{align}
Here we understand that $ \sup \emptyset = - \infty $. 
Note that $ T $ is finite almost surely. 
This random time $ T $ plays a role of stopping time in the sense that 
\begin{align}
\{ T=k \} \in \sigma(\phi_0,\phi_{-1},\ldots,\phi_k) 
\quad \text{for $ 0 \ge k+p-1 $}. 
\label{}
\end{align}
Since $ \sigma_p \circ \cdots \circ \sigma_1(V) $ is a singleton, 
we see that $ \phi_0 \circ \phi_{-1} \circ \cdots \circ \phi_{T} $ maps $ V $ into a singleton. 
Thus it holds that 
\begin{align}
X_0 = \phi_0 \circ \phi_{-1} \circ \cdots \circ \phi_{T}(x) 
\quad \text{a.s.} 
\label{eq: Xn=phin}
\end{align}
for all $ x \in V $. 
This shows the folloiwng: 
We pick a sequence $ f_0,f_{-1},\ldots $ from the law $ \mu $ 
up to the latest time $ T $ 
when $ (f_{T+p-1},\ldots,f_T) = (\sigma_p,\ldots,\sigma_1) $. 
Then the resulting site $ f_0 \circ f_{-1} \circ \cdots \circ f_T(x) $, 
which does not depend on the choice of $ x \in V $, 
is a sample point from the stationary law, 
which is as desired. 

This method can be applied to simulate a Gibbs distribution. 
In this case, a sync mapping law can be constructed 
with the help of monotonicity structure of the state space $ V $.

\begin{Rem}
The identity \eqref{eq: Xn=phin} 
implies that, for each $ n \in \bZ $, the random variable $ X_n $ 
is measurable with respect to $ \sigma(\phi_j:j \le n) $. 
One can ask what happens when $ \mu $ is not sync. 
The following theorem answers this question. 

\begin{Thm}[Yano \cite{Yrcp}] \label{thm: Yrcp}
Suppose that the Markov chain $ (X_n)_{n \in \bZ} $ is ergodic. 
and that $ \mu $ is not sync. 
Then, for each $ n \in \bZ $, the random variable $ X_n $ 
is not measurable with respect to $ \sigma(\phi_j:j \le n) $. 
\end{Thm}

For the proof of Theorem \ref{thm: Yrcp}, see \cite{Yrcp}. 
\end{Rem}

\section{Random walk and road coloring} \seclabel{sec: rc} 

Let us explain how our $ \mu $-random walk is related to road coloring. 

First, let us introduce some notations in graph theory. 
A finite directed graph is the pair $ (V,A) $ of finite sets $ V $ and $ A $ 
associated with mappings $ i:A \to V $ and $ t:A \to V $. 
Each element of $ V $ will be called a {\em site} (or a {\em node}) 
and each element $ a $ of $ A $ will be called a {\em (oneway) road} (or an {\em arrow}) 
which runs from $ i(a) $ to $ t(a) $. 
For $ a \in A $, the site $ i(a) $ (resp. $ t(a) $) will be called 
the {\em initial} (resp. {\em terminal}) site of $ a $. 
For $ x \in V $, the number of roads running from $ x $, namely, 
\begin{align}
O(x) = \sharp \cbra{ a \in A : i(a)=x } , 
\label{}
\end{align}
will be called the {\em outdegree} at the site $ x $. 
If $ O(x) $ does not depend on $ x \in V $, 
the directed graph $ (V,A) $ is called {\em of constant outdegree}. 
A {\em path} from $ x \in V $ to $ y \in V $ 
is a word $ w=(a_1,\ldots,a_n) $ of roads such that 
$ a_1 $ runs from $ x $ to $ i(a_2) $, 
$ a_2 $ to $ i(a_3) $, \ldots, $ a_{n-1} $ to $ i(a_n) $, 
and $ a_n $ to $ y $. 
The number $ L(w)=n $ is called the {\em length} of the path $ w=(a_1,\ldots,a_n) $. 
The directed graph $ (V,A) $ is called {\em strongly connected} 
if, for any $ x,y \in V $, there exists a path from $ x $ to $ y $. 
The directed graph $ (V,A) $ is called {\em aperiodic} 
if, for any $ x \in V $, the greatest common divisor 
of the set of $ L(w) \ge 1 $ among all paths $ w $ from $ x $ to itself is one. 

Second, we introduce some notations in road coloring. 
Suppose that $ (V,A) $ is of constant outdegree 
and denote the common outdegree by $ d $. 
A {\em road coloring} of $ (V,A) $ is a partition of $ A $ 
into $ d $ disjoint subsets $ C=\{ c^{(1)},\ldots,c^{(d)} \} $ 
such that, for each $ x \in V $, 
each {\em color} $ c^{(k)} $ contains one and only one road whose initial site is $ x $. 
For a finite sequence $ s=(c_1,\ldots,c_p) $ of elements of $ C $, 
a path $ w=(a_1,\ldots,a_p) $ is said to be {\em along $ s $} if 
$ a_k \in c_k $ for all $ k=0,1,\ldots,p $. 
The following notion originates Adler, Goodwyn and Weiss \cite{MR0437715}. 

\begin{Def} \label{def: sync C}
A road coloring $ C $ of $ (V,A) $ is called {\em sync} if 
there exists a finite sequence $ c_1,\ldots,c_p $ of elements of $ C $ 
such that all paths along $ (c_1,\ldots,c_p) $ have common terminal site. 
\end{Def}

Let us give an example. 

\begin{Ex} \label{ex: example2}
Let $ V=\{ 1,2,3 \} $ and $ A=\{ a^{(x,k)}:x \in V , \ k=1,2 \} $ 
and define the initial and terminal sites of each road as follows: 
\begin{align}
\text{
\begin{tabular}{c|cccccc}
$ a $    & $ a^{(1,1)} $ & $ a^{(2,1)} $ & $ a^{(3,1)} $ & $ a^{(1,2)} $ & $ a^{(2,2)} $ & $ a^{(3,2)} $ 
\\ \hline
$ i(a) $ & 1             & 2             & 3             & 1             & 2             & 3 
\\
$ t(a) $ & 3             & 3             & 1             & 3             & 1             & 2 
\end{tabular}
}
\label{eq: ex road}
\end{align}
Take the road coloring $ C = \{ c^{(1)},c^{(2)} \} $ defined by 
\begin{align}
c^{(1)} = \{ a^{(x,1)}:x \in V \} 
, \quad 
c^{(2)} = \{ a^{(x,2)}:x \in V \} . 
\label{eq: ex rc}
\end{align}
Now it is obvious that the road coloring $ C $ is sync; 
in fact, all paths along $ (c^{(1)},c^{(2)},c^{(2)},c^{(1)}) $ 
have common terminal site $ 3 $. 
\end{Ex}

Third, we recall the {\em road coloring problem}. 
If a directed graph $ (V,A) $ of constant outdegree admits a sync road coloring, 
then it is necessarily strongly connected and aperiodic. 
The converse was posed as a conjecture 
by Adler, Goodwyn and Weiss \cite{MR0437715}, 
which had been called the road coloring problem 
until it was completely solved by Trahtman \cite{MR2534238}. 

\begin{Thm}[Trahtman \cite{MR2534238}]
A directed graph which is of constant outdegree, 
strongly connected, and aperiodic, 
does admit a sync road coloring. 
\end{Thm}

Fourth, let us explain how to understand our $ \mu $-random walk by means of road coloring. 
Let $ \mu $ be a probability law on $ \Sigma $. 
Since $ \Sigma $ is a finite set, 
the support of $ \mu $ may be written as $ \{ \sigma^{(1)},\ldots,\sigma^{(d)} \} $. 
We define the set $ A $ of roads as the totality of 
$ a^{(x,k)} $ for $ x \in V $ and $ k=1,\ldots,d $ 
where $ a^{(x,k)} $ runs from $ x $ to $ \sigma^{(k)}(x) $. 
Thus the law $ \mu $ induces naturally 
the road coloring $ C = \{ c^{(1)},\ldots,c^{(d)} \} $ such that 
\begin{align}
c^{(k)} = \{ a^{(x,k)} : x \in V \} . 
\label{}
\end{align}
It is now obvious that 
the probability law $ \mu $ is sync in the sense of Definition \ref{def: sync mu}
if and only if the road coloring $ C $ is sync in the sense of Definition \ref{def: sync C}. 

For a $ \mu $-random walk $ (X,\phi) $, 
the process $ X $ moves from site to site in the directed graph $ (V,A) $ 
via the equation $ X_n=\phi_n(X_{n-1}) $, 
being driven by the colors of roads indicated by $ \phi $ 
which are randomly chosen from the road coloring $ C $ induced by $ \mu $. 
Thus we may call $ (X,\phi) $ a $ \mu $-random walk 
in the directed graph $ (V,A) $ subject to the road coloring $ C $. 

Let $ Y $ be a Markov chain 
and suppose that $ Y $ is realized as $ X $ of a $ \mu $-random walk $ (X,\phi) $ 
in the directed graph $ (X,\phi) $ subject to the road coloring $ C $ induced by $ \mu $. 
Then, to each edge $ (x,y) \in E(Y) $, 
there corresponds at least one road $ a $ which runs from $ x $ to $ y $. 
For example, consider Example \ref{ex: example} with $ p=0 $. 
In this case, we have $ \Supp(\mu) = \{ \sigma^{(3)},\sigma^{(4)} \} $, and 
\begin{align}
E(Y) = \{ (1,3),(2,1),(2,3),(3,1),(3,2) \} . 
\label{}
\end{align}
The set $ A $ of roads induced by $ \mu $ is given as $ A = \{ a^{(x,k)}:x \in V , k=1,2 \} $, 
where the initial and terminal sites of each road are given as \eqref{eq: ex road}. 
Then we find that the road coloring induced by $ \mu $ is 
nothing else but $ C = \{ c^{(1)},c^{(2)} \} $ 
given as \eqref{eq: ex rc} in Example \ref{ex: example2}, 
where we note that $ \sigma^{(3)} $ and $ \sigma^{(4)} $ 
correspond to $ c^{(1)} $ and $ c^{(2)} $, respectively. 
Note that there exist two roads $ a^{(1,1)} $ and $ a^{(1,2)} $ which run from 1 to 3, 
which are colored differently from each other. 
See Figure 1 below for the illustration. 

\begin{center}
\unitlength 0.1in
\begin{picture}( 20.0000, 16.0000)(  0.0000,-16.0000)
%
\special{pn 8}%
\special{ar 1000 200 200 200  0.0000000 6.2831853}%
%
\special{pn 8}%
\special{ar 200 1400 200 200  0.0000000 6.2831853}%
%
\special{pn 8}%
\special{ar 1800 1400 200 200  0.0000000 6.2831853}%
%
\special{pn 24}%
\special{pa 1170 320}%
\special{pa 1780 1200}%
\special{fp}%
\special{pa 1780 1200}%
\special{pa 1650 1150}%
\special{fp}%
\special{pa 1780 1200}%
\special{pa 1790 1060}%
\special{fp}%
\special{pa 390 1480}%
\special{pa 1610 1480}%
\special{fp}%
\special{pa 1610 1480}%
\special{pa 1500 1400}%
\special{fp}%
\special{pa 1610 1480}%
\special{pa 1500 1560}%
\special{fp}%
\special{ar 710 1310 1230 1230  5.1144264 6.2321348}%
%
\special{pa 1320 140}%
\special{pa 1190 180}%
\special{fp}%
\special{pa 1190 180}%
\special{pa 1280 310}%
\special{fp}%
%

\special{pn 8}%
\special{pa 1630 1260}%
\special{pa 1040 400}%
\special{fp}%
\special{pa 1630 1260}%
\special{pa 1500 1210}%
\special{fp}%
\special{pa 1630 1260}%
\special{pa 1640 1120}%
\special{fp}%
\special{pa 320 1220}%
\special{pa 880 370}%
\special{fp}%
\special{pa 880 370}%
\special{pa 740 420}%
\special{fp}%
\special{pa 880 370}%
\special{pa 900 510}%
\special{fp}%
\special{pa 1610 1320}%
\special{pa 390 1320}%
\special{fp}%
\special{pa 390 1320}%
\special{pa 500 1240}%
\special{fp}%
\special{pa 390 1320}%
\special{pa 500 1400}%
\special{fp}%
%

\put(10.0000,-2.2000){\makebox(0,0){1}}%
\put(1.9000,-14.1000){\makebox(0,0){2}}%
\put(18.0000,-14.0000){\makebox(0,0){3}}%
\put(5.2000,-6.7000){\makebox(0,0){$ a^{(2,2)} $}}%
\put(15.8000,-6.5000){\makebox(0,0){$ a^{(1,1)} $}}%
\put(12.5000,-9.0000){\makebox(0,0){$ a^{(1,2)} $}}%
\put(10.0000,-12.2000){\makebox(0,0){$ a^{(3,2)} $}}%
\put(10.0000,-15.8000){\makebox(0,0){$ a^{(2,1)} $}}%
\put(15.8000,-2.5000){\makebox(0,0){$ a^{(3,1)} $}}%
\end{picture}%
\\[3mm] Figure 1. 
\end{center}

\section{Conclusion} \seclabel{sec: conc} 

We have introduced a random walk in a finite set 
as a stochastic evolutionary process driven by an IID sequence of mappings. 
It can be understood as a random walk in a finite directed graph 
moving according to random road colors. 
Any ergodic Markov chain is proved to be realized, in a constructive way, 
by a random walk associated with a sync mapping law. 
The redundancy in random walk realization with a sync mapping law 
can be as close as desired to the minimum redundancy.

\label{lastpage-01}

\end{document}